# Active Learning at Scale: Investigating the Benefits of Peer Instruction in Undergraduate Mathematics

Raymond Vozzo[1]*, Stuart Johnson[1], Jonathan Tuke[1], and Tanya Evans[2]

[1] School of Computer and Mathematical Sciences, University of Adelaide, Adelaide, Australia; ORCID 0009-0003-7307-5396; raymond.vozzo@adelaide.edu.au; stuart.johnson@adelaide.edu.au; ORCID 0000-0002-1688-8951; simon.tuke@adelaide.edu.au

[2] Department of Mathematics, University of Auckland, Auckland CBD, Auckland, New Zealand; ORCID 0000-0001-5126-432X; t.evans@auckland.ac.nz

* **Corresponding authors:** raymond.vozzo@adelaide.edu.au

## Abstract

Active learning strategies have been widely recognised for their effectiveness in tertiary education, yet their implementation at scale, particularly in large first-year mathematics courses, presents considerable challenges. A common method for actively engaging students in large classes is through online quizzes, which may include structured peer instruction. In this study, we investigate the effect of having students answer quiz-style questions during class both with and without peer discussion in a first-year large mathematics course ($N$ = 550). We also investigate the short- and long-term effects of each protocol. Our findings indicate that peer instruction enhances student performance in mathematics in the following ways: First, when the responses to questions were measured before and after peer instruction the proportion of questions answered correctly increased by 0.2; second, when correct responses were compared to similar questions the following week the proportion correct increased by 0.34 (compared to 0.07 for the control); finally, when measured at the end of the semester the proportion of questions answered correctly increased by 0.42 (compared to 0.2 for the control). These results align with and build upon previous research on the benefits of peer instruction, indicating a potential long-term association with student knowledge retention of mathematical concepts. This suggests that peer discussion may do more than momentarily improve accuracy—it could also contribute to more durable learning.

## Introduction

Active learning has been proposed as a method for improving learning and outcomes in STEM areas (see for example: Freeman et al. (2014)). For large classes, traditional active learning techniques, such as group discussions or problem-solving exercises, often face logistical challenges. To overcome these, technology-enhanced methods have gained popularity, with one of the most well-known approaches being the use of "clickers" or audience response systems (Duncan & Mazur, 2005). These systems allow students to respond to quiz-style questions in real-time during lectures, providing instant feedback and enabling instructors to gauge student understanding on the spot.

One powerful variation of this approach is the integration of peer instruction, a technique first popularised by Eric Mazur at Harvard University (Crouch & Mazur, 2001). Peer instruction involves students initially answering a question on their own, followed by a peer discussion where students explain their reasoning and challenge each other's answers, before re-answering the same question. This process of discussion and reconsideration is thought to enhance comprehension, as students not only reinforce their own understanding but also benefit from exposure to alternative problem-solving approaches (Smith et al., 2009). The interaction among peers can potentially help clarify misconceptions and deepen students' grasp of the material.

Following the pioneering work of Mazur, who first employed peer instruction in physics, this approach has since been successfully adapted to a wide range of disciplines, including mathematics (Gao et al., 2025; Tullis & Goldstone, 2020). However, in mathematics education, only a small body of research has explored the use of peer instruction, mostly contrasting its effectiveness in comparison to traditional teaching methods, such as lectures. For example, studies have investigated the use of peer instruction in small classes in linear algebra (Teixeira, 2023) and calculus (Lucas, 2009; Pilzer, 2001), demonstrating its potential to enhance student motivation and reflection.

In this article, we explore the effectiveness of peer instruction in undergraduate mathematics courses, specifically within the contexts of large linear algebra and calculus classes. We compare student performance in answering questions during class, both with and without structured peer discussions. This design allows us to measure the added value of incorporating peer interactions into quiz-based active learning strategies. Our study explores the immediate impact of peer discussion on student performance and the extent to which these gains persist over time, providing insights into the delayed and long-term effectiveness of peer instruction in large mathematics courses.

The significance of this research lies in its potential to inform best practices for integrating active learning in large-scale mathematics classrooms. As many science and engineering faculties rely on first-year service courses to introduce foundational mathematical concepts to hundreds, or even thousands, of students, finding scalable methods that promote durable learning is essential. Peer instruction, combined with technology-driven tools, presents a promising solution to this challenge. By fostering student engagement and encouraging peer-to-peer collaboration, these methods can enhance the learning experience in large classes, ultimately contributing to better student outcomes in mathematics education.

## Theoretical Foundations

Peer instruction is a pedagogical strategy that has gained significant traction in postsecondary and tertiary education due to its demonstrated ability to enhance student engagement and conceptual understanding (Tullis & Goldstone, 2020). The approach, developed by Mazur (1997) in the context of university physics education, involves students first attempting a challenging problem individually before discussing their reasoning with peers and then reattempting the problem. A growing body of research has examined the cognitive and pedagogical mechanisms that drive the benefits of peer instruction.

# Cognitive and Educational Psychology Foundations

The effectiveness of peer instruction is underpinned by several well-established theories in cognitive and educational psychology, as well as core cognitive mechanisms known to support learning. In particular, the self-explanation effect, retrieval practice, conceptual change, and social development learning theory offer compelling explanations for how peer instruction can facilitate deeper understanding and promote long-term knowledge retention. These frameworks highlight the importance of active engagement, verbalisation of reasoning, and socially mediated learning in the development of robust conceptual understanding. The following sections examine each of these mechanisms in detail, outlining their relevance to peer instruction in mathematics education.

## Self-Explanation Effect

One of the most influential cognitive mechanisms supporting peer instruction is the self-explanation effect, first identified by Chi et al. (1994). Self-explanation occurs when learners articulate their reasoning while solving problems, leading to a deeper understanding by prompting them to infer missing knowledge, recognise inconsistencies, and refine their conceptual frameworks. The act of verbalising thought processes forces students to engage actively with the material, reinforcing their understanding. Peer instruction naturally incorporates self-explanation, as students must justify their reasoning to one another, often leading to the resolution of misconceptions and improved comprehension.

Research has shown that students who engage in self-explanation during problem-solving consistently outperform those who rely on passive study methods, such as simply reviewing worked examples (Renkl, 1997, 1999). In university mathematics education, this effect is particularly pronounced, as proof comprehension often requires advanced reasoning and logical deduction (Hodds et al., 2014). Self-explanation plays a crucial role in supporting knowledge integration and generalisation, which in turn improves understanding and future performance (Rittle-Johnson, 2024). When students explain their reasoning, they go beyond simply recalling mathematical concepts—they integrate and apply them in new contexts, facilitating the transfer of knowledge. There is strong empirical evidence that prompting learners to explain during mathematical learning leads to greater conceptual knowledge, procedural knowledge, and immediate procedural transfer. However, the evidence for long-term benefits of procedural transfer after a delay is limited (Rittle-Johnson, 2024). Scaffolding high-quality explanations through structured responses, designing prompts that focus on critical content, and ensuring learners explain both correct and incorrect information, further enhances the learning benefits of self-explanation (Berthold et al., 2009).

In peer instruction, the benefits of self-explanation could extend beyond individual cognitive gains. When engaging in peer instruction, students are often confronted with alternative perspectives that can challenge their initial understanding and lead to more robust knowledge construction. For example, students who initially hold incorrect conceptions may revise their thinking when forced to articulate and defend their answers, leading to greater accuracy (Tullis & Goldstone, 2020). This dynamic interaction is thought to reinforce the value of self-explanation within the context of peer instruction.

## Retrieval Practice and Spacing Effects

Retrieval practice, the process of actively recalling learned information, is one of the most well-documented cognitive strategies for enhancing memory and long-term retention. Research in cognitive psychology has repeatedly demonstrated that retrieving information strengthens memory traces, making future recall more efficient and reducing the likelihood of forgetting (Roediger & Butler, 2011; Roediger III & Karpicke, 2006). This effect is particularly relevant in mathematics, where proficiency depends on the ability to recall and apply foundational concepts fluently (Fuchs et al., 2021).

Peer instruction facilitates retrieval practice in multiple ways. First, when students attempt to answer a problem before discussion, they are required to access and apply their existing knowledge, reinforcing their understanding through self-testing. The subsequent discussion phase allows students to revisit and refine their responses, which further strengthens retention. Importantly, retrieval practice is most effective when it involves some level of difficulty—known as desirable difficulty—which encourages deeper cognitive engagement and long-term learning (Bjork & Bjork, 2011).

Additionally, peer instruction incorporates elements of spaced retrieval, which is the practice of recalling information at increasing intervals over time (Cepeda et al., 2009). Peer instruction is typically embedded within a unit of study, such as a semester-long course, and occurs on a weekly—or even more frequent—basis. Each mathematical concept is revisited multiple times throughout the semester through quizzes and peer discussions, ensuring that students repeatedly retrieve and apply their knowledge at spaced intervals. This structured repetition is thought to strengthen memory consolidation and deepen understanding over time (Pyc & Rawson, 2010). Research has shown that spaced retrieval enhances long-term retention far more effectively than massed practice, where information is reviewed in a single, concentrated session (Cepeda et al., 2006). By engaging in peer instruction regularly, students are not only encouraged to articulate and refine their understanding but are also prompted to retrieve and apply key concepts at multiple points throughout the course. This repeated, structured engagement makes spaced retrieval an inherent feature of peer instruction that should reinforce learning and support durable retention (Cepeda et al., 2006; Karpicke & Blunt, 2011).

Furthermore, discussions during peer instruction in groups often require students to explain their reasoning multiple times in response to different students' perspectives, effectively creating opportunities for repeated retrieval in varied contexts. This type of retrieval—retrieval variability—has been shown to improve transferability of knowledge, as students become adept at recognising when and how to apply their understanding in different problem scenarios (Butowska-Buczyńska et al., 2024).

In sum, the well-documented effects of spaced retrieval practice provide a strong theoretical foundation for understanding why peer instruction can enhance learning in mathematics education. By engaging in repeated retrieval of mathematical concepts through structured peer discussions, students can reinforce memory traces and strengthen their ability to recall and apply knowledge over time.

## Social Development Theory

Vygotsky's (1978) theory of social development provides a key theoretical underpinning for peer instruction. A central concept in this theory is the Zone of Proximal Development (ZPD), which describes the gap between what a learner can accomplish independently and what they can achieve with guidance from someone more knowledgeable. Learning is most effective within this zone, where support—known as scaffolding—can be provided to help the learner progress. Scaffolding refers to the temporary assistance or guidance offered to a learner to help them perform a task they could not complete alone; this support is gradually reduced as the learner gains confidence and independence. Peer instruction aligns with this framework by encouraging students to collaborate, enabling more knowledgeable peers to scaffold learning for others. Through such social interaction, learners co-construct understanding and advance their cognitive development.

This theoretical perspective is reinforced by empirical research from Chi (2013) that demonstrates that collaborative learning not only facilitates the sharing of explicit information between pairs of students but also promotes the co-construction of inferred knowledge. The findings indicate that interaction between learners can lead to the repair of mental models and the generation of new insights, reinforcing the value of peer dialogue in deepening conceptual understanding.

## Peer Discussion and Conceptual Change

Research on conceptual change in education suggests that one of the primary ways students overcome misconceptions is through cognitive conflict—when they encounter ideas that challenge their preexisting beliefs (Vosniadou, 2013a, 2013b). Peer instruction naturally fosters such cognitive conflict, as students are exposed to alternative perspectives during discussion. This aligns with the findings of Tullis and Goldstone (2020), who reported that initially incorrect students were more likely to revise their answers after discussion, particularly when confronted with explanations from confident peers. In general, this process of resolving cognitive conflict has been shown to be a powerful mechanism for conceptual change and deeper understanding (Chi, 2013).

However, empirical evidence also suggests that peer instruction can have unintended negative consequences in certain contexts. For instance, a recent study by Ospanbekov et al. (2024) in a university physics course found that for counterintuitive questions peer discussion did not lead to improved outcomes. In fact, a notable number of students changed their initially correct answers to incorrect ones after group discussion, resulting in a net decline in overall accuracy. This highlights the importance of considering question type and group dynamics when implementing peer instruction, as the mechanism of cognitive conflict may not always lead to productive resolution.

## Empirical Research on Peer Instruction in Mathematics Education

The body of empirical research investigating the impact of peer instruction in university mathematics courses remains relatively limited. Drawing on its success in physics education, Pilzer (2001) explored the integration of peer instruction into a small calculus course. After two semesters, he reported outcomes comparable to those observed in his physics courses, including significant improvements in students' reasoning skills and notable gains in knowledge retention. The evidence for these claims comes from a three-hour final review class, during which students revisited all peer instruction questions posed throughout the semester; for each question, over ninety percent of students selected the correct answer. Beyond suggested cognitive gains, Pitzer also suggested that peer instruction contributed to increased enthusiasm for calculus and a marked growth in students' confidence—outcomes that echoed his earlier experiences with peer instruction in physics contexts.

Similarly, Lucas (2009) implemented peer instruction using clickers in a small calculus course and noted perceived improvements in student engagement and conceptual understanding. However, his observations also highlighted a potential drawback: in some cases, students labeled by the author as "high-status"—typically those perceived as more knowledgeable—dominated discussions, sometimes persuading their peers to abandon correct answers in favor of incorrect ones. Analysis of classroom video recordings indicated that confident students often triggered passivity in others, who opted to listen rather than actively participate in discussions. To address this, Lucas recommended requiring students to first write down their individual answers before engaging in peer discussions, a strategy that appeared to mitigate one-sided interactions and foster more balanced exchanges.

More recently, Teixeira (2023) investigated the effects of peer instruction within a small linear algebra course ($N$ = 75), employing a novel approach that combined seminars with structured peer discussions. The study reported enhanced student understanding of abstract concepts and improved application skills, along with increased motivation and knowledge retention relative to traditional instruction. However, since peer instruction was implemented alongside several other instructional innovations, the causal impact of peer instruction alone remains uncertain.

An additional strand of research in mathematics education has focused on the design of peer instruction questions, aiming to optimize learning environment and outcomes. The research has focused on factors such as question clarity, cognitive load, and the alignment of questions with learning objectives (Cline et al.,

2013; Miller et al., 2006; Yamaoka et al., 2020). These studies highlight the importance of carefully crafting peer instruction activities to ensure they could promote engagement and support learning.

In summary, research on the effectiveness of peer instruction in mathematics education is limited, with no studies employing properly controlled conditions to support causal claims. While some evidence points to potential benefits, these findings should be interpreted cautiously. To advance the field, more nuanced studies employing randomised designs are needed to isolate the specific mechanisms through which peer instruction influences learning, while minimising confounding variables. This brings us to a critical open question that we outline in the following section.

## Open Question in Peer Instruction Research

An important consideration in peer instruction research is the extent to which learning gains persist over time. Tullis and Goldstone (2020) found that students who initially answered a question correctly tended to retain their answers more reliably than those who initially responded incorrectly. While confidence played a role in this trend, it did not fully account for the observed patterns, suggesting that peer instruction reinforces correct understandings while providing a pathway for conceptual change among students with initial misconceptions.

A critical question remains regarding the durability of these learning gains: do students who revise their answers during peer instruction retain their new, corrected understanding in the long term? Addressing this specific question is essential for determining the effectiveness of peer instruction in general. Our study aims to explore this by tracking student performance on similar problems over extended periods, providing insights into the durability of learning gains.

## Research Objective

Theoretical and empirical work suggests that peer instruction holds promise as an effective pedagogical approach in university mathematics education. Cognitive mechanisms such as self-explanation, retrieval practice, and socially mediated learning offer plausible explanations for how peer discussion might enhance student learning. While prior studies have indicated that peer instruction can improve immediate performance and conceptual reasoning, questions remain about the extent to which these gains are retained over time. Our study seeks to explore this dimension by examining student performance on related mathematical problems administered immediately, one week later, and at the end of the semester. Rather than aiming to establish causal claims, our goal is to provide preliminary, classroom-based evidence on the potential knowledge retention effects of peer instruction over multiple time points. In doing so, we aim to contribute to a more nuanced understanding of how structured peer interaction may support learning in mathematics, and to inform the ongoing development of evidence-based teaching practices in tertiary settings.

# Methodology

## Background

This study was conducted at a large, research-intensive Australian university in a first-year undergraduate mathematics course containing parallel streams of calculus and linear algebra, with a diverse cohort including mathematics, computer science, engineering and science students. Enrolments for the semester were approximately 550 students. For most students this is their second semester of university mathematics. The course utilises a flipped classroom model, with videos and notes on the course material available online in the learning management system; a one hour workshop each week for the entire class (in a very large workshop theatre) where the students participate in quizzes, with follow up discussion from the workshop where required; and weekly tutorials, where students work in groups of 4 or 5 at whiteboards solving problems. There are several components of assessment throughout the semester (including weekly written and online assignments, and an invigilated test) and a final written exam.

In the workshop, peer instruction is often used, where a question is asked, and students provide answers (using their phones or other devices to access Mentimeter) and are then invited to discuss the question with their peers and potentially change answers.

For the study, all students present in the workshop were assigned the same questions at the same time, but each question asked was randomly allocated one of two treatments:

- **control**, in which the question was posed and students have only one opportunity to answer, after which the solution would be revealed and explained if necessary. No particular attempt was made either to force students to discuss with their peers or to not discuss at all, they were allowed to answer in a natural way;

- **discussion**, in which students were instructed to give an initial answer without any interactions with their peers, with the collective answers shown to the class before giving the students a chance to discuss their answers, and change them if they wish. Then the solutions are shown.

The questions asked were mostly conceptual in nature, highlighting fundamental aspects of the material or particular topics that students typically have difficulty with. Generally it was expected that students could answer questions without the need for a great deal of calculation. Some examples of questions used can be found in Appendix.

The effect of peer instruction was measured in a number of different ways. In cases where students are asked to discuss questions and then given the opportunity to change their answer, we can compare the initial responses with those obtained after the discussion. To obtain a measure of the effect of peer instruction against a control group we measured responses over consecutive weeks. In each weekly class, students are first asked new questions about the same concepts as the questions from the previous week. These are designed to test conceptual understanding rather than recall of answers from the previous week (see Appendix for examples). This allows us to measure the effect of peer instruction by comparing performance on these questions for each of the "control" and "discussion" treatments.

In addition, some questions were repeated at the end of the semester in a revision class, allowing longer term retention of knowledge to be measured.

## Statistical Methodology

Each week, four new questions were asked in the quiz. These consisted of two algebra questions, and two calculus questions. Within each area (calculus or algebra), one question was a control question, and one was a discussion question. The allocation to either control or discussion was done by randomisation.

In weeks 2-11 four additional questions were asked, these were related to the questions from the previous week. In week 12, there were 15 questions that were all related to questions from previous weeks.

The form of the experiment is based on a matched pairs design, in that we have a pair of questions that were discussed in the previous week, and also a pair of questions that were not discussed and acted as the control. All of the analysis was performed in R (R Core Team 2023) using the program RStudio[1]. We fitted four models to address the effect of repetition on the knowledge retention rates of students. In each case, we fitted a generalised linear mixed-effects model (GLMM) with the outcome variable being the number of correct answers, and the predictors being the offering and whether the question had been reviewed—where appropriate. As well, to account for the repeated measurements, we included a random intercept for each question. The GLMM were fitted using the lme4 package (Bates et al., 2015). For each model, we then predicted the probability of getting a question correct as given in Table 1, Table 2, Table 3 and Table 4.

---
[1] https://posit.co/download/rstudio-desktop

All data collected was non-identifiable and the research was assessed as negligible risk in accordance with the National Statement on Ethical Conduct in Human Research.

## Results

Figure 1 gives an overview of the experimental design and indicates when questions were repeated. Some questions are the same question (brown triangle), while some are different questions, but the same concept (blue circle). The offering on the $x$-axis gives the week and the order, so for example 02-3 is the third question in the second week. The ID identifies each question. The first part gives the week, the letter indicates the subject (A is Algebra, while C is Calculus) and the number indicates whether there is more than one question in a subject, so W02A3 is the third question in Algebra for Week 2. We see that some questions were repeated once, for example W01A1, while some may be repeated up to four times (e.g. W04C2).

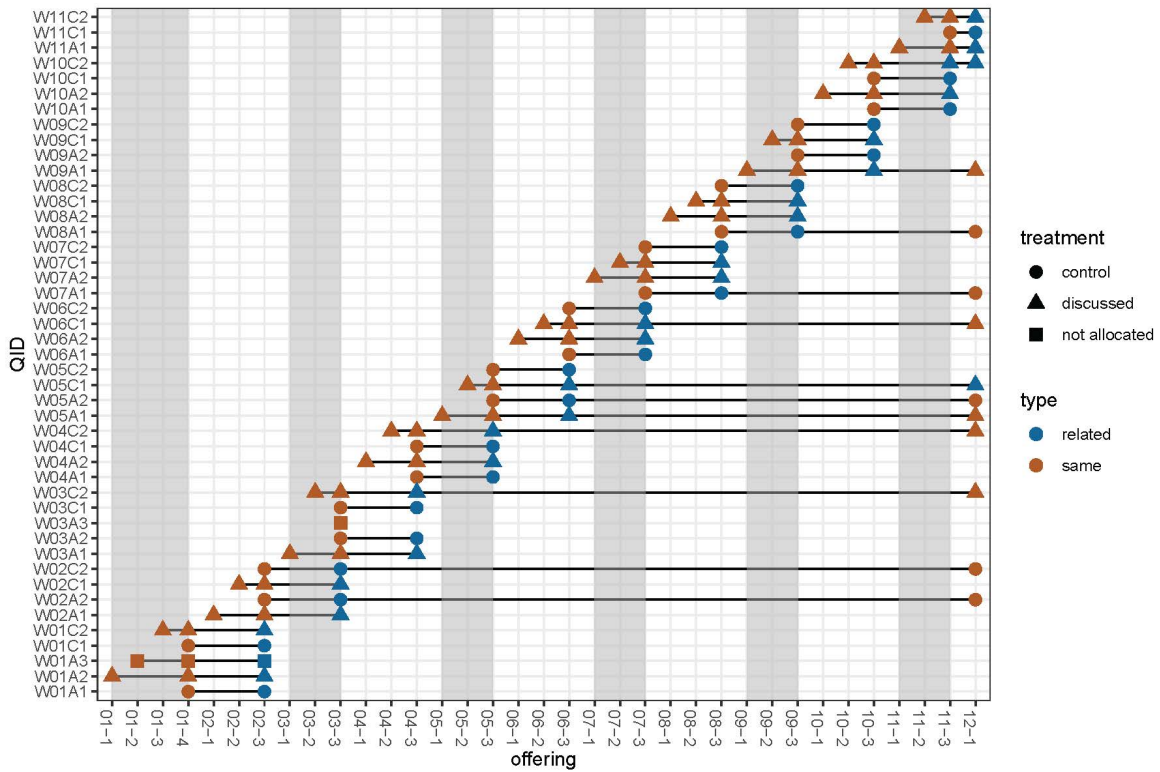

Figure 1: Figure of experimental design showing when questions were repeated. Some questions are the same question (brown triangle), while some are different questions, but the same concept (blue circle). The offering on the $x$-axis gives the week and the order; the Question ID identifies each question.

Figure 2 shows the proportion correct for each attempt of a question. The thick lines show the average proportion.

In more detail, we fitted four separate models:

- Repeat Model, which compares the first and second attempt at discussion questions in the same workshop.

- Short-term Model A, which compares the first attempt at control and discussion questions with the attempt at the related questions in the following week. Note that it is the first attempt at the discussion question being used.

- Short-term Model B, which compares the control and the second attempt at the discussion questions with the attempt at the related questions in the following week. Note that it is the second attempt at the discussion question being used.

- Long-term Model, which compares the first attempt at control and discussion questions with the attempt in Week 12.

The results of each of these models is shown in Table 1, Table 2, Table 3 and Table 4, respectively.

First, we see that over time there is an improvement in understanding as seen by an increased proportion of correct answers for second attempts and beyond. This is also seen in Table 1, Table 2, Table 3 and Table 4, where in all cases there is a positive improvement in proportion correct. Also, we see that the increase in understanding from the first attempt is much larger for discussion questions compared to the control questions (Table 4). This may be that we can see that the proportion correct for discussion questions on the first attempt is much lower than the control questions. Since the questions were randomly assigned to either the control or discussion group, this result is unexpected. Possible explanations for this outcome are explored in the Discussion section.

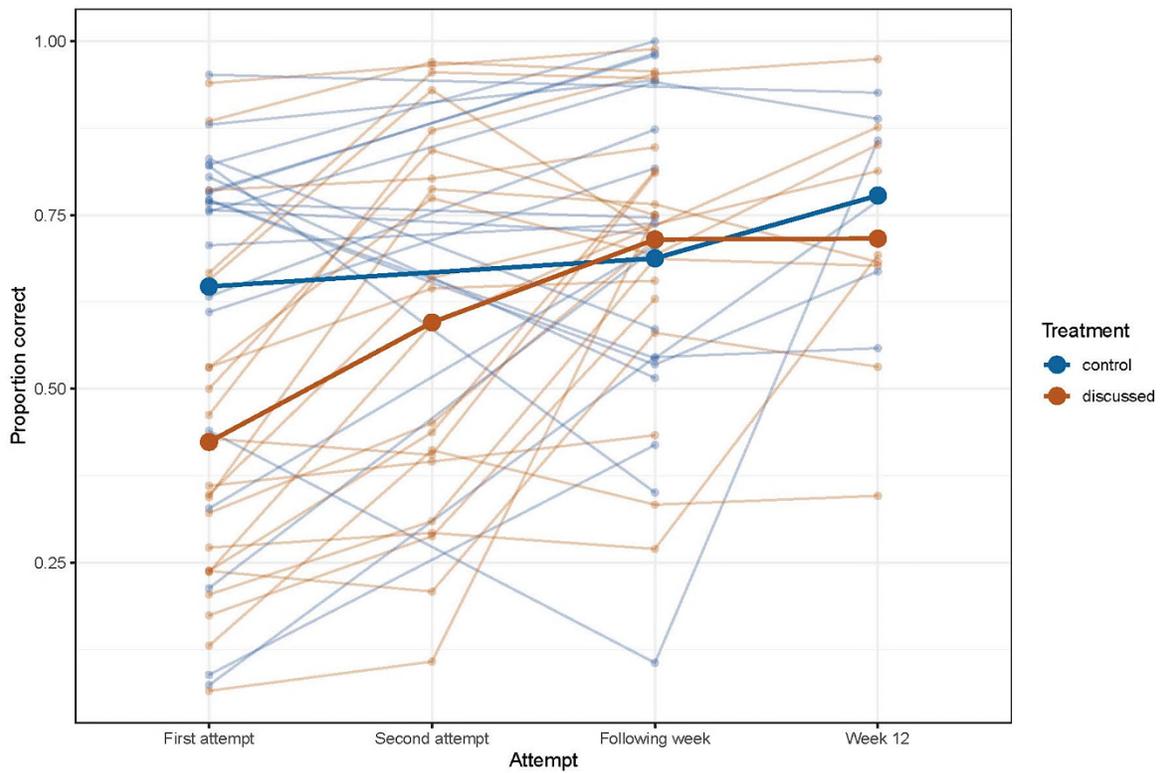

Figure 2: Proportion answered correctly for each attempt of the questions. The colour indicates if the question is a discussion or a control question. The lines correspond to a single question.

| Treatment | 1st | 2nd | Improvement |
|---|---|---|---|
| Discussed | 0.42 | 0.63 | 0.2 |

Table 1: Repeat Model. Proportion correct for questions at the first and then second attempt in the same workshop.

| Treatment | First week | Second week | Improvement |
|---|---|---|---|
| Control | 0.67 | 0.74 | 0.07 |
| Discussion | 0.43 | 0.77 | 0.34 |

Table 2: Short-term Model A. Proportion correct for questions at the first and then second attempt in the next workshop. Note that we are using the first attempt for the discussed questions in the first week.

| Treatment | First week | Second week | Improvement |
|---|---|---|---|
| Control | 0.67 | 0.74 | 0.07 |
| Discussion | 0.64 | 0.79 | 0.15 |

Table 3: Short-term Model B. Proportion correct for questions at the first and then second attempt in the next workshop. Note that we are using the second attempt for the discussed questions in the first week.

| Treatment | First offering | Week 12 | Improvement |
|---|---|---|---|
| Control | 0.55 | 0.75 | 0.2 |
| Discussion | 0.28 | 0.70 | 0.42 |

Table 4: Long-term Model. Proportion correct for questions at the first and then second attempt in the Week 12 workshop.

Table 5 contains the estimates and P-value for each of the models that we considered. We see a statistically significant effect (5% level) of offering and treatment (control versus discussion). Also, we see that there is a statistically significant interaction of offering and treatment for the Short-term and Long-term models. As the interpretation of two-way interactions are difficult to interpret, we advise looking at the predicted probabilities of getting the answers correct as given in Table 1, Table 2, Table 3 and Table 4.

| Model | Predictor | Estimate | P-value |
|---|---|---|---|
| Repeat | Repeat | 0.83 | $1.8 \times 10^{-33}$ |
| Short-term Model A | Discussion | −0.99 | $3.01 \times 10^{-3}$ |
|  | Second week | 0.34 | $1.05 \times 10^{-6}$ |
|  | Interaction | 1.16 | $1.11 \times 10^{-29}$ |
| Short-term Model B | Discussion | −0.13 | $7.25 \times 10^{-1}$ |
|  | Second week | 0.35 | $1.04 \times 10^{-6}$ |
|  | Interaction | 0.42 | $4.57 \times 10^{-5}$ |
| Long-term | Control | −1.13 | $2.34 \times 10^{-2}$ |
|  | Week 12 | 0.91 | $5.97 \times 10^{-11}$ |
|  | Interaction | 0.88 | $7.84 \times 10^{-6}$ |

Table 5: Coefficients and P-values for each of the three models fitted to look at the effect of offering and discussion on the probability of getting a question correct.

## Discussion

Our study adds to the growing body of evidence suggesting that peer instruction can support student learning in mathematics, aligning with and extending prior findings in the literature. Specifically, our results demonstrate that, at the group level, peer instruction significantly enhances immediate performance, as well as short-term (one week) and long-term (end-of-semester) retention of mathematical concepts. These effects can be understood more deeply through the theoretical lenses outlined in Theoretical Foundations: self-explanation effect, spaced retrieval effect, conceptual change theory, and social development theory.

First, we observed an immediate benefit from peer instruction, with the proportion of correct responses increasing by 0.2 after students engaged in discussion. This aligns with prior research, including Tullis and Goldstone (2020), which found that students who initially answered a question correctly tended to retain their answers more reliably than those who initially responded incorrectly. While student confidence played a role in this trend, it did not fully account for the observed patterns, suggesting that peer instruction provides students with opportunities to clarify misunderstandings and refine their reasoning in real time—processes that are central to the conceptual change framework. Verbalising one's thinking to peers not only exposes inconsistencies or misconceptions but also facilitates self-explanation, a cognitive process that encourages learners to make explicit connections between concepts, thus deepening their understanding

(Rittle-Johnson, 2024). The dialogic nature of peer instruction also reflects principles from social development theory (Vygotsky, 1978), where interaction with more knowledgeable peers or co-constructing knowledge within the zone of proximal development fosters cognitive growth.

More importantly, our study extends previous research by providing evidence for the persistence of these learning gains over time. One of the key open questions raised by Tullis and Goldstone (2020) was whether errors corrected during peer instruction remain corrected on later assessments. While we did not track individual students' responses across timepoints, our group-level results provide preliminary insight. When students were tested on similar questions the following week, the overall proportion of correct responses increased by 0.34 for questions that had been discussed during peer instruction, compared to only 0.07 for the control group. Although we cannot rule out the possibility that some students who initially answered correctly later changed to incorrect responses, this net improvement indicates that the mechanisms activated during immediate peer discussion—namely self-explanation, conceptual restructuring, and socially mediated learning—continue to reinforce and consolidate knowledge over time. This interconnected process suggests that peer instruction not only enhances initial comprehension but also fosters durable retention of mathematical concepts.

This claim is also supported by our end-of-semester results that showed a greater increase in the proportion of correct responses for the peer instruction questions (0.42) compared to the control questions (0.2), suggesting that the benefits of peer discussion may extend across a longer time frame. While this pattern is encouraging, it is important to note that, as with the short-term findings, we did not track individual performance over time. Therefore, the observed group-level improvement cannot confirm whether the same students retained or deepened their understanding. Nonetheless, the relative difference in performance suggests that peer instruction might contribute to more persistent learning, potentially by fostering deeper cognitive engagement with underlying concepts.

This possibility aligns with the notion of conceptual change, where learners reorganise or replace pre-existing knowledge structures through meaningful engagement with ideas. Peer instruction may provide a context in which such restructuring is more likely to occur—particularly when students are prompted to explain, justify, and evaluate mathematical reasoning with others. These socially mediated interactions, as described in social development theory, may help learners articulate and refine their thinking in ways that support both immediate comprehension and longer-term conceptual development.

Taken together, our findings provide preliminary evidence that peer instruction may enable durable mathematical learning. However, given the nature of our study design and the classroom context, these conclusions must be interpreted with caution. Understanding the impact of peer instruction on student learning is inherently complex, particularly in natural classroom settings where experimental control is limited. Our study was not a randomised controlled experiment in which students were assigned to either a peer instruction or a non-discussion condition. As a result, we cannot draw strong causal conclusions. Ethical and logistical constraints in real classroom environments made such a design infeasible, as withholding the opportunity for discussion from some students would raise pedagogical concerns.

Additional limitations also stem from the variability in how students engaged with the peer instruction process. The workshop theatre layout meant that students could choose whether or not to sit with peers and engage in discussion. Some may have discussed the questions deeply, while others may have passively observed the class response patterns displayed on Mentimeter without verbal interaction. Consequently, it is possible that some answer changes were influenced more by conformity to the group than by genuine conceptual shifts.

Moreover, the distinction between "control" and "peer instruction" conditions was not always clean-cut. While the peer instruction questions involved structured discussion and explicit prompts from lecturer to first answer individually, the control and repeated questions may not have been entirely free of peer interaction. Students were not prevented from discussing during control questions, and in a naturalistic learning environment, occasional informal discussion is difficult to eliminate. Thus, rather than a strict dichotomy, these conditions reflect different levels and structures of peer engagement—more accurately conceptualised as "structured peer instruction" vs "business-as-usual".

This recognition led to the development of Short-term Model B (Table 3), which compares first attempts on control questions with second (post-discussion) attempts on peer instruction questions. The model suggests that students benefitted from first attempting a question independently before engaging in discussion, a finding aligned with well-documented effects of self-explanation and retrieval practice. However, future studies could improve internal validity by ensuring that review questions administered the following week are completed without discussion. This would provide clearer evidence of the retention effects observed after a one-week delay.

Further work is also needed to explore how different physical classroom environments might influence the dynamics and effectiveness of peer instruction. Smaller, flat-floor rooms may better support collaborative dialogue than large lecture-style theatres. Future iterations of this research might compare the impact of peer instruction across varied spatial contexts to determine whether the setting moderates the effectiveness of peer interaction.

Finally, it would be valuable to examine whether the benefits of peer instruction differ based on individual factors such as prior knowledge, confidence, or engagement style. Students with stronger conceptual foundations may approach peer discussions differently than those with less background knowledge, potentially resulting in different learning outcomes. Understanding these nuances could support more targeted and equitable instructional strategies and deepen our understanding of how peer instruction operates across diverse learner profiles.

## Acknowledgments

The authors would like to thank Nickolas Falkner for helpful feedback on the interpretation of the results.

## Appendix: Examples of Questions

We provide some examples of questions used during the study that demonstrate the sort of thing that students were typically asked. The examples below show both conceptual questions and computational questions. Note that if any computations were required, they were both fairly simple and designed to test understanding of a particular technique or result. In each case two versions of a question are presented, the initial one, and the related question which was asked the following week. The related questions are designed to test the same concept, and to require genuine understanding of the concept, rather than being able to be correctly answered purely by knowing the answer to the original question.

**Example 1** Which of the following are linear transformations?
- A. $F(x,y) = (x,y)$
- B. $F(x,y) = (0,0)$
- C. $F(x,y) = (1,1)$
- D. $F(x,y) = (xy, 0)$
- E. $F(x,y) = (y,y)$.

**Related question** Which of the following are linear transformations?
- A. $F(x, y) = (y, x)$
- B. $F(x, y) = (0, 1)$
- C. $F(x, y) = (x + y, 0)$
- D. $F(x, y) = (xy, xy)$
- E. $F(x, y) = (y, y + 1)$.

**Example 2** Suppose $A = \begin{bmatrix} 2 & * \\ * & 3 \end{bmatrix}$ is not invertible. What are the eigenvalues of $A$?
- A. 0
- B. 2
- C. 3
- D. 5
- E. 6
- F. There is not enough information to determine this

**Related question** Suppose $A = \begin{bmatrix} -1 & * \\ * & 2 \end{bmatrix}$ is not invertible. What are the eigenvalues of $A$?
- A. $-2$
- B. $-1$
- C. 0
- D. 1
- E. 3
- F. There is not enough information to determine this

**Example 3** Fill in the blank.
$\sum_{n=1}^{\infty} a_n$ converges _______ $\lim_{n \to \infty} a_n = 0$.
- A. $\Rightarrow$
- B. $\Leftarrow$
- C. $\Leftrightarrow$

**Related question** Fill in the blank.
Suppose that $a_n > 0$ is a decreasing sequence. Then $\sum_{n=1}^{\infty} (-1)^n a_n$ converges _______ $\lim_{n \to \infty} a_n = 0$.
- A. $\Rightarrow$
- B. $\Leftarrow$
- C. $\Leftrightarrow$

**Example 4** Consider the graph of the function $f(x, y)$ below.

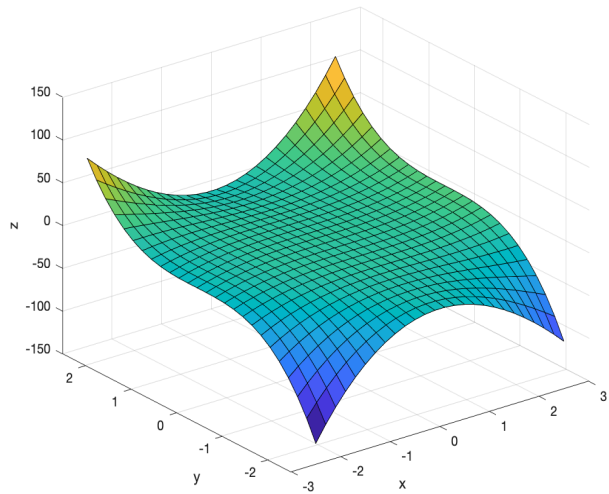

Which of the following is $f_x$ and which is $f_y$?

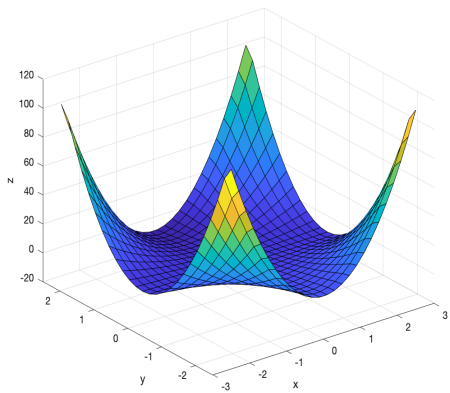

A.

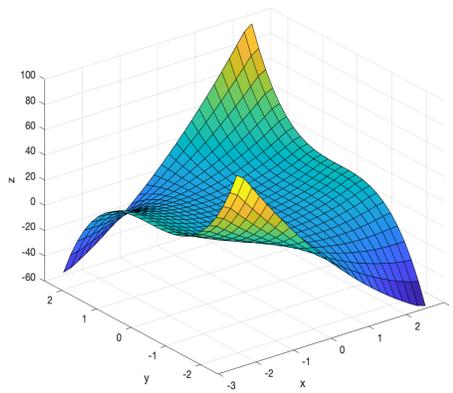

B.

**Related question** Consider the graph of the function $f(x, y)$ below.

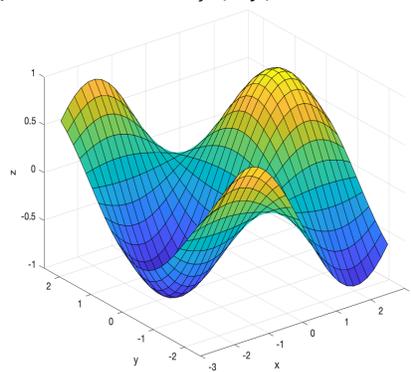

Which of the following is $f_x$ and which is $f_y$?

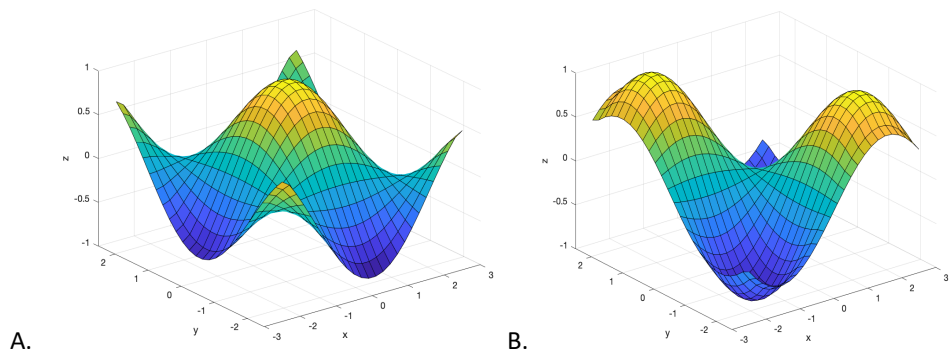

A.          B.

## Data Availability Statement

The data that support the findings of this study are openly available in figshare at http://doi.org/10.25909/28852628.